# On Orbits of the Ring $Z_n^m$ under the Action of the Group $SL(m,Z_n)$

Petr Novotný, Jiří Hrivnák

Department of Physics,
Faculty of Nuclear Sciences and Physical Engineering, Czech
Technical University, Břehová 7, 115 19 Prague 1, Czech
Republic

Ing. Petr Novotný     fujtajflik@seznam.cz     222311333
Ing. Jiří Hrivnák     hrivnak@post.cz          222311333

**Keywords**: Ring, finite group

## Abstract
We consider the action of the finite matrix group $SL(m,Z_n)$ on the ring $Z_n^m$. We determine orbits of this action for $n$ arbitrary natural number. It is a generalization of the task which was studied by A.A. Kirillov for $m=2$ and $n$ prime number.

## 1. Introduction

The important role of symmetries in classical and quantum physics is well known. We focus on so called discrete quantum physics; it means that the corresponding Hilbert space is finite dimensional [1,2]. Well known are also 2×2 Pauli matrices, besides spanning real Lie algebra su(2), they form a fine grading of sl(2,C). The fine gradings of a given Lie algebra are preferred bases which yield quantum observables with additive quantum numbers.

The generalized $n \times n$ Pauli matrices were described in [3]. For $n=3$ these 3×3 Pauli matrices form one of four non-equivalent gradings of sl(3,C). Other fine gradings are Cartan decomposition and the grading which corresponds to Gell-Mann matrices [4,5]. The symmetries of the fine grading of sl($n$,C) associated with these generalized Pauli matrices were studied only recently in [6]. This work pointed out the importance of the finite group $SL(2,Z_n)$ as the group of symmetry of the Pauli gradings. The additive quantum numbers, mentioned above, form in this case the finite associative additive ring $Z_n \times Z_n$. The action of $SL(2,Z_n)$ on $Z_n \times Z_n$ then represents the symmetry transformations of Pauli gradings of sl($n$,C). The orbits of this action form such points in $Z_n \times Z_n$ which can be reached by symmetries.

For the purpose of so called graded contractions [7], it became convenient to study the action of $SL(2,Z_n)$ on various types of Cartesian products of $Z_n$ [8]. Note that the orbits of $SL(2,Z_p)$ on $Z_p^2$, where $p$ is prime number were considered in [9] §16.3. The purpose of this article is to generalize this result to orbits of $SL(m,Z_n)$ on $Z_n^m$ where $m,n$ are arbitrary natural numbers.

## 2. Action of the group $SL(m,Z_n)$

In the whole article we shall use the following notation: $N:=\{1,2,3,\ldots\}$ denotes the set of all natural numbers and $P:=\{2,3,5,\ldots\}$ denotes the set of all prime numbers. Let $n$ be a natural number, then the set $\{0,1,\ldots,n-1\}$ forms, together with operations $+_{\mathrm{mod}\,n}$, $\cdot_{\mathrm{mod}\,n}$, an associative commutative ring with unity. We will denote this ring, as usual, by $Z_n$. It is well known that for $n$ prime the ring $Z_n$ is a field.

Let us consider $m,n$ be arbitrary natural numbers. We denote by
$$Z_n^m = \underbrace{Z_n \times Z_n \times \ldots \times Z_n}_{m}$$
the Cartesian product of $m$ rings $Z_n$. It is clear that $Z_n^m$ with operations $+_{\mod n}$, $\cdot_{\mod n}$ defined elementwise is an associative commutative ring with unity again. It contains divisors of zero and we call its elements **row vectors** or **points**. Furthermore we call the zero element $(0,\ldots,0)$ **zero vector** and denote it simply by 0.

We denote by $Z_n^{m,m}$ the set of all $m \times m$ matrices with elements in the ring $Z_n$. For $k \in N$ and $A \in Z_n^{m,m}$ we will denote by $(A)_{\mod k}$ a matrix which arose from matrix A after application of operation modulo $k$ on its elements.

In the following we shall frequently use a product on the set $Z_n^{m,m}$ defined as matrix multiplication together with operation modulo $n$, i.e.
$$A, B \in Z_n^{m,m} \to (AB)_{\mod n} . \tag{2.1}$$
This product is, due to associativity of matrix multiplication, associative again and the set $Z_n^{m,m}$ equipped with this product forms a semigroup. If we take matrices $A, B \in Z_n^{m,m}$, such that $\det(A) = \det(B) = 1 \pmod{n}$, then $\det((AB)_{\mod n}) = 1 \pmod{n}$ holds. It follows that the subset of $Z_n^{m,m}$ formed by all matrices with determinant equal to unity modulo $n$ is a semigroup.

**Definition 2.1:** For $m, n \in N$, $n \geq 2$ we define $SL(m, Z_n) := \{A \in Z_n^{m,m} \mid \det A = 1 \pmod{n}\}$.

Now we show that $SL(m, Z_n)$ with operation (2.1) forms a group. Because $SL(m, Z_n)$ is a semigroup, it is sufficient to show that there exists a unit element and a right inverse element. Unit matrix is clearly the unit element. In order to find a right inverse element consider the following equation
$$AA^{adj} = \det(A)I. \tag{2.2}$$
The symbol $A^{adj}$ denotes the adjoint matrix defined by $(A^{adj})_{i,j} := (-1)^{i+j} \det A(j,i)$, where $A(j,i)$ is the matrix obtained from matrix A by omitting the j-th row and the i-th column. The equation (2.2) holds for an arbitrary matrix, hence it holds for matrices from $SL(m, Z_n)$, and evidently holds after application of operation modulo $n$ on both sides. Consequently, for $A \in SL(m, Z_n)$, we have
$$AA^{adj} = I \pmod{n}, \quad \text{i.e.} \quad (AA^{adj})_{\mod n} = I .$$
Therefore $A^{adj}$ is the right inverse element corresponding to the matrix A, and consequently $SL(m, Z_n)$ is a group.

The group $SL(m, Z_n)$ is finite and its order was computed by You Hong and Gao You in [10] (see also [11], p. 86). If $n \in N$, $n \geq 2$ is written in the form $n = \prod_{i=1}^{r} p_i^{k_i}$, where $p_i$ are distinct primes, then according to [10], the order of $SL(m, Z_n)$ is

$$|SL(m, Z_n)| = n^{m^2-1} \prod_{i=1}^{r} \prod_{j=2}^{m} \left(1 - \frac{1}{p_i^j}\right). \tag{2.3}$$

Let $G$ be a group and $X \neq \emptyset$ a set. Recall that a mapping $\psi: G \times X \to X$ is called a **right action** of the group $G$ on the set X if the following conditions hold for all elements $x \in X$:
1. $\psi(gh, x) = \psi(g, \psi(h, x))$ for all $h, g \in G$.
2. $\psi(e, x) = x$, where $e$ is the unit element of $G$.

Let ψ be an action of a group *G* on a set X. A subset of *G*, {$g \in G \mid \psi(g,a)=a$} is called a **stability subgroup** of the element $a \in X$. A subset of X, {$b \in X \mid \exists g \in G, \ b = \psi(g,a)$} is called an **orbit** of the element $a \in X$ with respect to the action ψ of the group *G*.

Let us note that if ψ is an action of a group *G* on a set X then relation ~ defined by formula
$$a,b \in X, \ a \sim b \Leftrightarrow \exists g \in G, \ \psi(g,a)=b \qquad (2.4)$$
is an equivalence on the set X and the corresponding equivalence classes are orbits.

**Definition 2.2:** For $m,n \in \mathbb{N}$, $n \geq 2$ we define a right action ψ of the group $SL(m,Z_n)$ on the set $Z_n^m$ as right multiplication of the row vector $a \in Z_n^m$ by the matrix $A \in SL(m,Z_n)$ modulo *n*:
$$\psi(A,a) := (aA)_{\mathrm{mod}\ n}$$
Henceforth we will omit the symbol mod *n* and write this action simply as *a*A.

## 3. Orbits for *n=p* prime number

The purpose of this section is to describe orbits of the ring $Z_p^m$ under the action of the group $SL(m,Z_p)$, where *p* is prime. Trivially, for *m* =1 is $SL(1,Z_p)=\{(1)\}$ and any orbit has the form {*a*} for $a \in Z_p$. Consequently we will further consider $m \geq 2$. It is clear that the zero element can be transformed by the action of $SL(m,Z_p)$ to itself only, thus it forms a one-point orbit and its stability subgroup is the whole $SL(m,Z_p)$. Let us take a nonzero element for instance $(0,\ldots,0,1) \in Z_p^m$, and find its orbit. An arbitrary matrix A from $SL(m,Z_p)$ acts on this element as follows

$$(0,\ldots,0,1) \begin{pmatrix} A_{1,1} & A_{1,2} & \cdots & A_{1,m} \\ \vdots & \vdots & \vdots & \vdots \\ A_{m-1,1} & A_{m-1,1} & \cdots & A_{m-1,m} \\ A_{m,1} & A_{m,1} & \cdots & A_{m,m} \end{pmatrix} = (A_{m,1}, A_{m,2}, \ldots, A_{m,m}) \quad (\mathrm{mod}\ p).$$

Thus the orbit of element $(0,\ldots,0,1)$ contains last row of any matrix from $SL(m,Z_p)$. It follows from det(A) = 1 that these rows cannot be zero and we show that they can be equal to an arbitrary nonzero element from $Z_p^m$. Let $(A_{m,1},A_{m,2},\ldots,A_{m,m}) \in Z_p^m$ be a nonzero element, that means $\exists j \in \{1,2,\ldots,m\}$ such that $A_{m,j} \neq 0$, then the matrix A can be chosen with determinant equal to 1. Without loss of generality consider *j* = 1:

$$A = \begin{pmatrix} 0 & & & \\ \vdots & & B & \\ 0 & & & \\ A_{m,1} & A_{m,2} & \cdots & A_{m,m} \end{pmatrix}, \text{ where } B = \mathrm{diag}(1,\ldots,1,(-1)^{1+m}(A_{m,1})^{-1}).$$

Here $(A_{m,1})^{-1}$ denotes the inverse element to $A_{m,1}$ in the field $Z_p$.
We conclude that in the case of *n=p* prime there are only two orbits:

1. one-point orbit represented by the zero element $(0,\ldots,0,0)$

2. $(p^m-1)$-point orbit $Z_p^m \setminus \{0\}$ represented by the element $(0,\ldots,0,1)$

## 4. Orbits for *n* natural number

We consider an arbitrary natural number *n* of the form

$$n = \prod_{i=1}^{r} p_i^{k_i},$$

where $p_i$ are distinct primes and $k_i$ are natural numbers.

The action of the group $SL(m,Z_n)$ on the ring $Z_n^m$ was established in definition 2.2 as a right multiplication of a row vector from $Z_n^m$ by a matrix from $SL(m,Z_n)$ modulo *n*. We define an equivalence induced by this action on the ring $Z_n^m$ according to (2.4). Elements $a=(a_1,a_2,\ldots,a_m)$, $b=(b_1,b_2,\ldots,b_m) \in Z_n^m$ are equivalent $a \sim b$ if and only if there exists $A \in SL(m,Z_n)$ such that $aA=b$ i.e. $\sum_{j=1}^{m} a_j A_{i,j} = b_i \pmod{n}$, $\forall i \in \{1,2,\ldots,m\}$.   (4.1)

**Definition 4.1:** Let $\sim$ be the equivalence on $Z_n^m$ defined by (4.1). For any divisor *d* of *n*, we will denote by $\text{Or}_{m,n}(d)$ the class of equivalence (orbit) containing the point $(0,\ldots,0,(d)_{\text{mod } n})$, i.e.

$$\text{Or}_{m,n}(d) = \{a \in Z_n^m \mid a \sim (0,\ldots,0, (d)_{\text{mod } n})\}. \quad (4.2)$$

Note that the orbit $\text{Or}_{m,n}(n)$ contains only the zero vector, because the zero vector can be transformed by the action of $SL(m,Z_n)$ only to itself. We shall see later that any orbit in $Z_n^m$ has the form (4.2).

**Definition 4.2: A greatest common divisor** of the element $a = (a_1,a_2,\ldots,a_m) \in Z_n^m$ and the number $n \in N$ is the greatest common divisor of all components of the element *a* and the number *n* in the ring of integers Z. We denote it by

$$\gcd(a,n) := \gcd(a_1,a_2,\ldots,a_m,n). \quad (4.3)$$

**Lemma 4.3:** The action of the group $SL(m,Z_n)$ on the ring $Z_n^m$ preserves the greatest common divisor of an arbitrary element $a \in Z_n^m$ and the number *n*, i.e.

$$\gcd(aA,n) = \gcd(a,n) \qquad \forall a \in Z_n^m, \forall A \in SL(m,Z_n).$$

**Proof:** It follows from $aA = (\sum_{i=1}^{m} a_i A_{i,1}, \ldots, \sum_{i=1}^{m} a_i A_{i,m})$ and $\gcd(a,n) \mid \sum_{i=1}^{m} a_i A_{i,j}$, $\forall j \in \{1,2,\ldots,m\}$ that $\gcd(a,n) \mid \gcd(aA,n)$, i.e. the greatest common divisor cannot decrease during this action. If we take an element $aA$ and a matrix $A^{-1}$ we obtain $\gcd(aA,n) \mid \gcd(aAA^{-1},n) = \gcd(a,n)$ and together with first condition we have $\gcd(aA,n) = \gcd(a,n)$.

QED

**Corollary 4.4:** For any divisor *d* of *n* the orbit $\text{Or}_{m,n}(d)$ is a subset of $\{a \in Z_n^m \mid \gcd(a,n)=d\}$.

We will show that the orbit $\text{Or}_{m,n}(1)$ is equal to the set $\{a \in Z_n^m \mid \gcd(a,n)=1\}$. From corollary 4.4 we know that $\text{Or}_{m,n}(1)$ is the subset of $\{a \in Z_n^m \mid \gcd(a,n)=1\}$ and we prove that they have the same number of elements. At first we determine the number of points in $\text{Or}_{m,n}(1)$. For this purpose we determine the stability subgroup of the element $(0,\ldots,0,1)$. It is obviously formed by matrices of the form

$$A = \begin{pmatrix} A_{1,1} & A_{1,2} & \cdots & A_{1,m} \\ \vdots & \vdots & \vdots & \vdots \\ A_{m-1,1} & A_{m-1,1} & \cdots & A_{m-1,m} \\ 0 & 0 & \cdots & 1 \end{pmatrix}, \qquad \det(A) = 1 \pmod{n}.$$

Expansion of this determinant gives

$$1 = \det(A) = (-1)^{m+m} \det A(m,m) = \det A(m,m) \pmod{n}.$$

Therefore the stability subgroup of the point $(0,\ldots,0,1)$ is:

$$S := \left\{ A = \begin{pmatrix} & B & & A_{1,m} \\ & & & A_{2,m} \\ & & & \vdots \\ 0 & 0 & \cdots & 1 \end{pmatrix} \in SL(m,Z_n) \,|\, B \in SL(m-1,Z_n) \right\},$$

and its order is
$$|S| = n^{m^2-m-1} \prod_{i=1}^{r} \prod_{j=2}^{m-1} (1 - p_i^{-j}). \tag{4.4}$$

According to the Lagrange theorem, product of the order and the index of an arbitrary subgroup of a given finite group is equal to the order of this group. If we define on the group $SL(m,Z_n)$ a left equivalence induced by the stability subgroup $S$ by formula

$$A, B \in SL(m,Z_n) \qquad A \approx_S B \Leftrightarrow AB^{-1} \in S,$$

then we obtain equivalence classes of the form $SB = \{AB \mid A \in S\}$, $B \in SL(m,Z_n)$, i.e. right cosets from $SL(m,Z_n)/S$. The number of these cosets is, by definition, the index of the subgroup $S$. These cosets correspond one-to-one with points of the orbit which includes the point $(0,\ldots,0,1)$. Therefore the index of the stability subgroup $S$ is equal to the number of points in this orbit. Similar calculation can be done for an arbitrary point in an arbitrary orbit. Thus we have the following proposition.

**Proposition 4.5:** The number of elements in an orbit is equal to the order of the group $SL(m,Z_n)$ divided by the order of the stability subgroup of an arbitrary element in this orbit.

Using (2.3) and (4.4) we obtain that the number of points in the orbit $\mathrm{Or}_{m,n}(1)$ is equal to

$$|\mathrm{Or}_{m,n}(1)| = n^m \prod_{i=1}^{r} (1 - p_i^{-m}). \tag{4.5}$$

Now we will determine the number of all elements in $Z_n^m$, which have the greatest common divisor with the number $n$ equal to unity. This number is equal to the Jordan function.

**Definition 4.6:** For $m \in N$ a mapping $\varphi_m : N \to N$ defined by

$$\varphi_m(n) = |\{ a \in Z_n^m \mid \gcd(a,n) = 1 \}| \tag{4.6}$$

is called the **Jordan function** of the order $m$.

We present, without proof, some basic properties of the Jordan function which can be found in [12].

**Proposition 4.7:** For the Jordan function $\varphi_m$ of the order $m \in \mathbb{N}$ and for any $n \in \mathbb{N}$ holds:

1. $$\varphi_m(n) = n^m \prod_{p|n,\, p \in P} (1 - p^{-m}) \tag{4.7}$$

2. $$\sum_{d|n,\, d \in \mathbb{N}} \varphi_m(d) = n^m \tag{4.8}$$

3. $$\varphi_m\left(\frac{n}{d}\right) = |\{ a \in \mathbb{Z}_{\frac{n}{d}}^m \mid \gcd(a, \tfrac{n}{d}) = 1 \}| = |\{ a \in \mathbb{Z}_n^m \mid \gcd(a,n) = d \}|. \tag{4.9}$$

The number of all elements in $\mathbb{Z}_n^m$, which are co-prime with $n$, given by the first property of the Jordan function $\varphi_m(n)$ (4.7), is equal to the number of points in the orbit $\mathrm{Or}_{m,n}(1)$. Therefore the orbit $\mathrm{Or}_{m,n}(1)$ is formed by all elements in $\mathbb{Z}_n^m$ which are co-prime with $n$.

**Proposition 4.8:** For $m, n \in \mathbb{N}$, $m \geq 2$ holds $\mathrm{Or}_{m,n}(1) = \{ a \in \mathbb{Z}_n^m \mid \gcd(a,n) = 1 \}$.

## 4.1 Orbits for $n = p^k$ power of a prime

Let us now consider $n$ of the form $n = p^k$, where $p$ is a prime number and $k \in \mathbb{N}$, and determine orbits in this case.

**Definition 4.1.1:** For $j \in \mathbb{N}$, $j \leq k$, we define a mapping $F^j : \mathbb{Z}_{p^k}^m \to \mathbb{Z}_{p^k}^m$ by the formula
$$F^j(a) = (p^j.a)_{\bmod p^k} \qquad \text{for any } a \in \mathbb{Z}_{p^k}^m.$$

**Lemma 4.1.2:** Let $a$ and $b$ be two equivalent elements from $\mathbb{Z}_{p^k}^m$ and $j \leq k$. Then the elements $F^j(a)$ and $F^j(b)$ are equivalent as well.

**Proof:** Let $a, b \in \mathbb{Z}_{p^k}^m$, $a \sim b$. It follows from definition of equivalence $\sim$ that there exists a matrix $A \in SL(m, \mathbb{Z}_{p^k})$ such that $aA = b$. Consequently $F^j(aA) = F^j(b)$, where
$$F^j(aA) = (p^j aA)_{\bmod p^k} = (p^j a)_{\bmod p^k} (A)_{\bmod p^k} = F^j(a)A.$$
Since we have $F^j(a)A = F^j(b)$ and therefore $F^j(a) \sim F^j(b)$. \hfill QED

**Proposition 4.1.3:** Any orbit in the ring $\mathbb{Z}_{p^k}^m$ has the form
$$\mathrm{Or}_{m,p^k}(p^j) = \{ a \in \mathbb{Z}_{p^k}^m \mid \gcd(a, p^k) = p^j \}, \qquad 0 \leq j \leq k,$$
and consists of $|\mathrm{Or}_{m,p^k}(p^j)| = \varphi_m(p^{k-j})$ points.

**Proof:** From Lemma 4.1.2 it is clear that $F^j$ maps the orbit $\mathrm{Or}_{m,p^k}(1)$ into the orbit $\mathrm{Or}_{m,p^k}(p^j)$ and from Corollary 4.4 we have
$$F^j(\mathrm{Or}_{m,p^k}(1)) \subset \mathrm{Or}_{m,p^k}(p^j) \subset \{ a \in \mathbb{Z}_{p^k}^m \mid \gcd(a, p^k) = p^j \}.$$

Conversely,
$$\{a\in Z_{p^k}^m \mid \gcd(a,p^k) = p^j\}=\{ p^j a \mid a\in Z_{p^{k-j}}^m, \gcd(a,p^{k-j})=1\} \subseteq$$
$$\subseteq \{ (p^j a)_{\bmod p^k} \mid a\in Z_{p^k}^m, \gcd(a,p^k) = 1\} = F^j(\text{Or}_{m,p^k}(1)).$$

Thus we have
$$F^j(\text{Or}_{m,p^k}(1))=\text{Or}_{m,p^k}(p^j) =\{a\in Z_{p^k}^m \mid \gcd(a,p^k) = p^j \}. \qquad \text{QED}$$

## 4.2 Orbits for $n=pq$, $\gcd(p,q)=1$

Let us now consider $n$ of the form $n=pq$, where $p,q\in N$ are co-prime numbers. In this case it will be very useful to apply the Chinese remainder theorem [13].

**Theorem 4.2.1:** (Chinese remainder theorem)
Let $a_1, a_2 \in Z$. Let $p_1, p_2 \in N$ be co-prime numbers. Then there exists $x\in Z$, such that
$$x=a_i \pmod{p_i}, \forall i=1,2.$$
If $x$ is a solution, then $y$ is a solution if and only if
$$x=y \pmod{p_1 p_2}.$$

**Definition 4.2.2:** For $p,q\in N$, $\gcd(p,q)=1$ we define a mapping G: $Z_{pq}^m \to Z_p^m \times Z_q^m$ by the formula
$$G(a) := ((a)_{\bmod p}, (a)_{\bmod q}) \quad \text{for any } a\in Z_{pq}^m,$$

and a mapping g: $SL(m,Z_{pq}) \to SL(m,Z_p) \times SL(m,Z_q)$ by the formula
$$g(A) := ((A)_{\bmod p}, (A)_{\bmod q}) \quad \text{for any } A\in SL(m,Z_{pq}).$$

It is clear from definition that G,g are homomorphisms and the Chinese remainder theorem implies that G,g are one-to-one correspondences. Thus we have the following proposition.

**Proposition 4.2.3:** The mapping G is an isomorphism of rings and the mapping g is an isomorphism of groups.

Further we determine orbits on the Cartesian product of rings $Z_p^m \times Z_q^m$. For this purpose we define action of the Cartesian product of groups $SL(m,Z_p) \times SL(m,Z_q)$ on ring $Z_p^m \times Z_q^m$ by the formula
$$aA = (a_1,a_2)(A_1,A_2)=((a_1 A_1)_{\bmod p}, (a_2 A_2)_{\bmod q})$$

for any $a=(a_1,a_2)\in Z_p^m \times Z_q^m$ and any $A=(A_1,A_2)\in SL(m,Z_p) \times SL(m,Z_q)$. It follows from definition of this action that orbits in $Z_p^m \times Z_q^m$ are Cartesian products of orbits in $Z_p^m$ and $Z_q^m$.

**Proposition 4.2.4:** Let $p,q\in N$ be co-prime numbers. Then the mapping G provides one-to-one correspondence between orbits in $Z_{pq}^m$ and Cartesian products of orbits in $Z_p^m$ and $Z_q^m$. Moreover, if $p_1|p$, $q_1|q$ and the orbits $\text{Or}_{m,p}(p_1)$, $\text{Or}_{m,q}(q_1)$ are of the form
$$\text{Or}_{m,p}(p_1)=\{a\in Z_p^m \mid \gcd(a,p)=p_1\}, \quad \text{Or}_{m,q}(q_1)=\{a\in Z_q^m \mid \gcd(a,q)=q_1\},$$
then
$$\text{Or}_{m,pq}(p_1 q_1) = G^{-1}(\text{Or}_{m,p}(p_1) \times \text{Or}_{m,q}(q_1)) = \{a\in Z_{pq}^m \mid \gcd(a,pq) = p_1 q_1\}.$$

**Proof:** At first we prove that G and $G^{-1}$ preserve equivalence, i.e.
$$a \sim b \Leftrightarrow G(a) \sim G(b) \quad \text{for all} \quad a,b \in Z_{pq}^m.$$
From definition of equivalence we have
$$a \sim b \Leftrightarrow \exists A \in SL(m, Z_{pq}), \ aA = b \Leftrightarrow G(aA) = G(b),$$
where
$$G(aA) = ((aA)_{\text{mod } p}, (aA)_{\text{mod } q}) = ((a)_{\text{mod } p}, (a)_{\text{mod } q})((A)_{\text{mod } p}, (A)_{\text{mod } q}) = G(a)g(A).$$
Because G and g are one-to-one correspondences we obtain
$$a \sim b \Leftrightarrow aA = b \Leftrightarrow G(a)g(A) = G(b) \Leftrightarrow G(a) \sim G(b).$$
Since the mapping G is an isomorphism and G, $G^{-1}$ preserve equivalence, orbits in the ring $Z_{pq}^m$ one-to-one corresponds to orbits in the ring $Z_p^m \times Z_q^m$, and these are Cartesian products of orbits on $Z_p^m$ and $Z_q^m$.

Now remain to prove that the orbit $Or_{m,pq}(p_1q_1)$ corresponds to the orbit $Or_{m,p}(p_1) \times Or_{m,q}(q_1)$. It follows from the Chinese remainder theorem that G maps the set $\{a \in Z_{pq}^m | \gcd(a,pq) = p_1q_1\}$ on the set $\{(a_1,a_2) \in Z_p^m \times Z_q^m | \gcd(a_1,p) = p_1, \gcd(a_2,q) = q_1\}$, which is equal to the orbit $Or_{m,p}(p_1) \times Or_{m,q}(q_1)$. Therefore the set $\{a \in Z_{pq}^m | \gcd(a,pq) = p_1q_1\}$ forms an orbit and from Corollary 4.4 follows that
$$Or_{m,pq}(p_1q_1) = \{a \in Z_{pq}^m | \gcd(a,pq) = p_1q_1\}. \qquad \text{QED}$$

As a corollary of Propositions 4.1.3 and 4.2.4 we obtain the following theorem.

**Theorem 4.9:** Consider the decomposition of the ring $Z_n^m$, $m \geq 2$ into orbits with respect to the action of the group $SL(m, Z_n)$. Then

i) any orbit is equal to the orbit $Or_{m,n}(d)$ for some divisor $d$ of $n$, i.e.
$$Z_n^m = \bigcup_{d|n} Or_{m,n}(d);$$

ii) $Or_{m,n}(d) = \{a \in Z_n^m | \gcd(a,n) = d \}$;

iii) the number of points $|Or_{m,n}(d)|$ in $d$-orbit is given by the Jordan function
$$|Or_{m,n}(d)| = \varphi_m\left(\frac{n}{d}\right) = \left(\frac{n}{d}\right)^m \prod_{pd|n, p \in P}(1 - p^{-m}).$$

## Conclusion:

We have stepwise determined the orbits on the ring $Z_n^m$ with respect to the action of the group $SL(m, Z_n)$. At first we have proceeded in the same way as Kirillov in [9] and we have obtained the orbits in the case of $n$ prime number. In this case there are only two orbits, the first is one-point orbit formed by zero element and the second one is formed by all nonzero elements. The next step was the case of $n = p^k$ power of prime, there we found $k+1$ orbits characterized by the greatest common divisor of their elements and number $n$. Finally the orbits for an arbitrary natural number $n$ were found. Our results are summarized in Theorem 4.9.

## Acknowledgements

We would like to thank prof. Jiří Tolar, prof. Miloslav Havlíček and doc. Edita Pelantová for numerous stimulating and inquisitive discussions.